\documentclass[]{asl}

\usepackage{enumerate}

\title[$\Sigma^0_{3}$ Determinacy and $\boldsymbol{\Pi}^1_2$ monotone induction]{$\Sigma^0_{3}$ Determinacy and $\boldsymbol{\Pi}^1_2$ monotone induction}

\author{Sherwood J. Hachtman}
\revauthor{Hachtman, Sherwood J.}

\address{Department of Mathematics, Statistics, and Computer Science\\
University of Illinois at Chicago\\
Chicago, IL 60613, USA}

\email{hachtma1@uic.edu}

\thanks{We thank Itay Neeman for encouragement during the preparation of this work, and Donald A.\ Martin for referring us to his paper \cite{MaIND}, which inspired the proofs in Section~\ref{sec:pmi}.  The results here appear as Chapter 2 of the author's Ph.D.\ dissertation.}

\newtheorem{Theorem}{Theorem}[section]

\newtheorem{Lemma}[Theorem]{Lemma}
\newtheorem{Corollary}[Theorem]{Corollary}
\newtheorem{Claim}[Theorem]{Claim}

\theoremstyle{definition}
\newtheorem{Definition}[Theorem]{Definition}

\theoremstyle{remark}
\newtheorem{Remark}[Theorem]{Remark}

\renewcommand{\P}{\mathcal{P}}

\newcommand{\om}{\omega}

\newcommand{\ZF}{\mathsf{ZF}}
\newcommand{\ACA}{\mathsf{ACA}}
\newcommand{\CA}{\textsf{-CA}}
\newcommand{\ATR}{\mathsf{ATR}}

\newcommand{\KP}{\mathsf{KP}}
\newcommand{\KPI}{\mathsf{KPI}}

\newcommand{\Gdet}{\Gamma\text{-DET}}
\newcommand{\fone}{f_{\text{I}}}
\newcommand{\ftwo}{f_{\text{II}}}
\newcommand{\Mi}{\mathcal{M}_{\operatorname{I}}}
\newcommand{\Mii}{\mathcal{M}_{\operatorname{II}}}
\newcommand{\wfo}{\operatorname{wfo}}
\newcommand{\wfp}{\operatorname{wfp}}
\newcommand{\ON}{\operatorname{ON}}

\newcommand{\plto}{\rightharpoonup}
\newcommand{\rst}{\upharpoonright}

\newcommand{\M}{\mathcal{M}}

\newcommand{\Th}{\operatorname{Th}}

\newcommand{\PPCA}{\boldsymbol{\Pi}^1_2\text{-}\mathsf{CA}_0}
\newcommand{\PCA}{\boldsymbol{\Pi}^1_1\text{-}\mathsf{CA}_0}
\newcommand{\DCA}{\boldsymbol{\Delta}^1_2\text{-}\mathsf{CA}_0}
\newcommand{\SSAC}{\boldsymbol{\Sigma}^1_2\text{-}\mathsf{AC}_0}
\newcommand{\phil}{\prec_\Phi}

\newcommand{\R}{\mathbb{R}}

\newcommand{\DET}{\operatorname{-DET}}

\newcommand{\bgP}{{\boldsymbol{\Pi}}}

\newcommand{\gS}{\Sigma}
\newcommand{\gP}{\Pi}
\newcommand{\gD}{\Delta}

\newcommand{\PMI}{\boldsymbol{\Pi}^1_2\text{-}\mathsf{MI}}
\newcommand{\GMI}{{\Gamma}\text{-}\mathsf{MI}}
\newcommand{\MI}{\text{-}\mathsf{MI}}

\newcommand{\la}{\langle}
\newcommand{\ra}{\rangle}
\newcommand{\nGS}{\neg \Game \boldsymbol{\Sigma}^0_3}

\newcommand{\ID}{\textsf{-ID}}
\newcommand{\TR}{\textsf{-TR}}
\newcommand{\TI}{\textsf{-TI}}
\begin{document}

\begin{abstract}Building on recent work of Philip Welch (\cite{We}, \cite{We2}), we prove that (lightface) $\Sigma^0_3$ determinacy is equivalent to the existence of a wellfounded model satisfying the axiom scheme of (boldface) $\boldsymbol{\Pi}^1_2$ monotone induction.\end{abstract}

\maketitle

\section{Introduction}
In this paper, we isolate the strength of $\Sigma^0_3$ determinacy in terms of a natural theory in second order arithmetic.  Namely, we show that $\gS^0_3\DET$ is equivalent over $\PCA$ to the existence of a countably-coded $\beta$-model of \emph{$\bgP^1_2$ monotone induction}.  

There is a great deal of precedent for equivalences between determinacy in low levels of the Borel hierarchy and axioms of inductive definition.  In one of the first studies in reverse mathematics, Steel \cite{St} proved over $\mathsf{RCA}_0$ that $\ATR_0$ is equivalent to both $\gD^0_1\DET$ and $\gS^0_1\DET$.  Tanaka \cite{TaDelta} showed over $\ACA_0$ that $\Delta^0_2\DET$ is equivalent to $\Pi^1_1\TR$, and  \cite{Ta} that over $\ATR_0$, $\gS^0_2\DET$ is equivalent to $\gS^1_1\MI$.  MedSalem and Tanaka \cite{MedTanDelthree} established equivalences over $\ATR_0$ between $k{-}\gP^0_2\DET$ and $[\gS^1_1]^k\ID$, an axiom allowing inductive definitions using combinations of $k$-many $\gS^1_1$ operators; furthermore, they showed over $\gP^1_3\TI$ that $\Delta^0_3\DET$ is equivalent to $[\gS^1_1]^{\text{TR}}\ID$, an axiom allowing inductive definition by combinations of transfinitely many $\gS^1_1$ operators.  Further results were given by Tanaka and Yoshii \cite{TYDetInd} characterizing the strength of determinacy for pointclasses refining the difference hierarchy on $\gP^0_2$, again in terms of axioms of inductive definition.

Just beyond these pointclasses we have $\gS^0_3$, where an exact characterization of strength has been elusive.  The sharpest published bounds on this strength were given by Welch \cite{We}, who showed that although $\gS^0_3\DET$ (and more) is provable in $\gP^1_3\CA_0$, $\gD^1_3\CA_0$ (even augmented by $\textsf{AQI}$, an axiom allowing definition by arithmetical quasi-induction) cannot prove $\gS^0_3\DET$.  On the other hand, Montalb\'{a}n and Shore \cite{MS} showed that $\gS^0_3\DET$ (and indeed, any true $\gS^1_4$ sentence) cannot prove $\gD^1_2\CA_0$.  This situation is further clarified by the same authors in \cite{MScon}, where they show (among other things) that $\gS^0_3\DET$ implies the existence of a $\beta$-model of $\gD^1_3\CA_0$.

Welch $\cite{We2}$ went on to give a characterization of the ordinal stage at which winning strategies in $\gS^0_3$ games are constructed in $L$. There, the least ordinal $\gamma$ so that every $\Sigma^0_3$ game is determined with a winning strategy definable over $L_{\gamma}$ is shown to be the least $\gamma$ for which there exists an illfounded admissible model $\M$ with an infinite descending sequence of nonstandard levels of $L$ that fully $\Sigma_2$-reflect to standard levels below $\gamma$, and so that $\wfo(\M)=\gamma$ (see Definition~\ref{IDN}).

Montalb\'{a}n asks, as Question 28 of \cite{reverseopen}, for a precise classification of the proof-theoretic strength of $\Sigma^0_3\DET$.  
In light of the work of Welch and Montalb\'{a}n-Shore, it appeared plausible that $\gS^0_3\DET$ is actually equivalent to the existence of a $\beta$-model of some natural theory in second order arithmetic.  We felt that the ordinal $\gamma$ appearing as the wellfounded ordinal of Welch's nonstandard structure should be characterized as the least so that $L_\gamma$ satisfies some form of monotone induction.  This is what we show: $L_\gamma$ is the minimal model closed under $\bgP^1_2$ monotone inductive definitions, and indeed, $\gS^0_3\DET$ is equivalent over $\PCA$ to the existence of such a model.

All published proofs of $\Sigma^0_3$ determinacy trace back to Morton Davis's \cite{Da} (for relevant definitions, see Section~\ref{sec:det}).  Let $A \subseteq \om^\om$ be a $\Sigma^0_3$ set, so that $A = \bigcup_{k \in \omega} B_k$ for some recursively presented sequence $\la B_k\ra_{k \in \om}$ of $\Pi^0_2$ sets.  The idea behind the proof that the game $G(A)$ is determined is a simple one: if Player I does not have a winning strategy, then Player II refines to a quasistrategy  $W_0$ so that no infinite plays in $W_0$ belong to $B_0$, and so that $W_0$ doesn't forfeit the game for Player II (in the sense that I has no winning strategy in $G(A;W_0)$).  Having done this, Player II plays inside $W_0$ and at all positions of length 1, refines further to a $W_1$ which avoids $B_1$ without forfeiting the game.  Then refine to $W_2$ at positions of length 2, and so on.  The ultimate refinement of the sequence $W_0, W_1, W_2\dots$ of quasistrategies gives a winning quasistrategy for Player II in $G(A)$, since every infinite play must eventually stay in each $W_n$, and so avoid each $B_n$.

The key claim that makes this proof work is Lemma~\ref{davislemmalow} below, which asserts that whenever Player I does not have a winning strategy in $G(A;T)$, then for all $k$, there is such a quasistrategy $W_k$ for II.  Welch's characterization amounts to an analysis of the way in which these $W_k$ first appear in $L$.  Namely, if $T \in L_{\gamma}$ is such that I doesn't have a winning strategy for $G(A,T)$ in $L_\gamma$, then the assumed reflection of the ordinals of $\M$ ensures there is a quasistrategy $W \in L_\gamma$ as in the conclusion of the lemma.  Furthermore, Welch defines a game that is won by the $\Pi^0_3$ player, but for which there can be no winning strategy in $L_\gamma$.  In this situation it is necessarily the case that the quasistrategies $W_k$ from which II's winning strategy is built are constructed cofinally in $L_\gamma$, and the common refinement is only definable over the model $L_\gamma$.

Welch's proof of determinacy is difficult, and the quasistrategies of interest are obtained in something of a nonstandard way.  Our present aim is to give a more constructive account of the way in which the quasistrategies $W_k$ arise.  There is a relatively straightforward way in which $W_k$ can obtained by iteration of a certain monotone operator.  The complexity of this operator is $\nGS$ in the parameter $T$ (here $\Game$ is the game quantifier as defined in \cite{DST} 6D; $\neg\Gamma$ denotes the dual pointclass of $\Gamma$).  It seemed natural to conjecture, then, that the ordinal $\gamma$ is in some sense a closure ordinal for these monotone inductive definitions (indeed, Welch makes the conjecture in \cite{We2} that $\gamma$ is $o(\Game \Pi^0_3)$, the closure ordinal of \emph{non}-monotone $\Game \Pi^0_3$ inductive definitions; it appears however, based on the results here, that $\gamma < o(\Game \Pi^0_3)$).

In what follows, we denote subsets of $\omega$ by capital Roman letters $X,Y,Z$, elements of $\omega$ by lowercase Roman letters from $i$ up to $n$, ordinals by lowercase Greek $\alpha, \beta, \gamma...$ and reals (elements of $\R = \om^\om$) by $w,x,y,z$.

\begin{Definition}\label{def:indops} Let $\Gamma$ be a pointclass.  $\GMI$ is the axiom scheme asserting, for each $\Phi:\P(\om) \to \P(\om)$ which is a $\Gamma$ operator, i.e.
 \[
  \{ \la n, X \ra \mid n \in \Phi(X) \} \in \Gamma,
 \]
that is monotone, i.e.
 \[
   (\forall X, Y) X \subseteq Y \to \Phi(X) \subseteq \Phi(Y),
 \]
that there exists an ordinal $o(\Phi)$ and sequence $\la \Phi^\xi \ra_{\xi \leq o(\Phi)}$ such that, setting $\Phi^{<\xi} = \bigcup_{\zeta < \xi} \Phi^{\zeta}$, we have \begin{itemize}
  \item for all $\xi \leq o(\Phi)$, 
  $   
    \Phi^\xi = \Phi(\Phi^{<\xi}) \cup \Phi^{<\xi},
  $
  \item $\Phi^{o(\Phi)} = \Phi^{<o(\Phi)}$, and
  \item $o(\Phi)$ is the least ordinal with this property.
\end{itemize}
$\Phi^{o(\Phi)}$ is the \emph{least fixed point of $\Phi$}, denoted $\Phi^\infty$.
\end{Definition}
There is a prewellorder $\phil$ with field $\Phi^\infty \subseteq \om$ naturally associated with the sequence $\la \Phi^\xi \ra_{\xi \leq o(\Phi)}$.  Namely, set $m \phil n$ if and only if the least $\xi$ with $m \in \Phi^\xi$ is less than the least $\zeta$ with $n \in\Phi^{\zeta}$.

We are interested in the case that $\Gamma$ is one of $\boldsymbol{\Pi}^1_2$, $\Pi^1_2(z)$ for a real $z$.  We regard $\Pi^1_2(z)\MI$ as being formalized in the language $\mathsf{L}_2$ of second order arithmetic: It is the schema asserting the existence of the prewellorder $\phil$, for each $\Pi^1_2(z)$ monotone operator $\Phi$.  Note for such $\Phi$, the relation ``$X = \phil$", as a relation holding of $X \in \P(\om \times \om)$, is arithmetical in $\Sigma^1_2(z)$.

It will follow from our analysis that $\gamma$ is least so that $L_\gamma \models (\nGS)\MI$.  However, we argue to more directly show something stronger: that $L_\gamma \models \PMI$.  (We remark that, as $\gS^0_3\DET$ fails in $L_\gamma$, $\nGS$ and $\Game \bgP^0_3$ do not coincide there.)   For a deep analysis of iterated monotone $\nGS$ and $\bgP^1_2$ operators, we encourage the reader to consult Martin's \cite{MaIND}.

In Section~\ref{sec:det}, we show that winning strategies in $\Sigma^0_3$ games are definable over any $\beta$-model of $\PMI$.  In Section~\ref{sec:pmi}, we prove that Welch's infinite depth $\Sigma_2$-nestings furnish us with such $\beta$-models.  We complete this circle of implications in Section~\ref{sec:wel} by reproducing Welch's lower bound argument in the base theory $\PCA$ to show that $\Sigma^0_3$ determinacy implies the existence of an infinite depth $\Sigma_2$-nesting.  We conclude with an analysis of the $\boldsymbol{\Pi}^1_2$ relations which are correctly computed in $L_\gamma$: these are precisely the relations $\nGS$ in parameters from $L_\gamma$.

\section{Proving determinacy}\label{sec:det}
Let $T$ be a (non-empty) tree with no terminal nodes; $[T]$ denotes the set of infinite branches of $T$, and for $p \in T$, $T_p$ denotes the subtree of $T$ with stem $p$, that is, $T_p = \{q \in T \mid q \subseteq p \vee p \subseteq q \}$.  For a set $A \subseteq [T]$, the \emph{game on $T$ with payoff $A$,} denoted $G(A;T)$, is defined as the infinite perfect information game in which two players, I and II, alternate choosing successive nodes of a branch $x$ of $T$; we call such an infinite branch a \emph{play}.  Player I wins if $x \in A$; otherwise, Player II wins.  We write $G(A)$ for $G(A;{\om^{<\om}})$.

A \emph{strategy for I} in a game on $T$ is a partial function $\sigma: T \plto X$ that assigns to an even-length position $s \in T$ a legal move $x$ for I at $s$, that is, $x \in X$ so that $s^\frown \langle x \rangle \in T$.  We require the domain of $\sigma$ to be closed under legal moves by II as well as moves by $\sigma$; note then that due to the presence of terminal nodes in the tree, I needn't have a strategy at all.  Strategies for II are defined analogously.  If an infinite play $x$ can be obtained by playing against a strategy $\sigma$, we say $x$ is \emph{according to} or {compatible with} $\sigma$.  We say a strategy $\sigma$ is \emph{winning} for Player I (II) in $G(A;T)$ if every play according to $\sigma$ belongs to $A$ ($[T] \setminus A$).  A game $G(A;T)$ is \emph{determined} if one of the players has a winning strategy.  For a pointclass $\Gamma$, $\Gdet$ denotes the statement that $G(A;\om^{<\om})$ is determined for all $A \subseteq \om^\om$ in $\Gamma$.  

We furthermore define a \emph{quasistrategy for Player II in $T$} to be a subtree $W \subseteq T$, again with no terminal nodes, that does not restrict Player I's moves, in the sense that whenever $p \in W$ has even length, then every 1-step extension $p^\frown \la s \ra \in T$ belongs to $W$.  A quasistrategy may then be thought of as a multi-valued strategy.  (Similar definitions of course can be made for Player I, but at no point will we need to refer to  quasistrategies for Player I.)

Quasistrategies are typically obtained in the following fashion: if Player I does not have a winning strategy in $G(A;T)$, then setting $W$ to be the collection of $p \in T$ so that I doesn't have a winning strategy in $G(A;T_p)$, we have that $W$ is a quasistrategy for II in $T$.  This $W$ is called \emph{II's nonlosing quasistrategy in $G(A;T)$}.

\begin{Theorem}\label{detfrommi} Let $\M$ be a $\beta$-model of $\PMI$.  Then for any $\Sigma^0_3(z)$ set $A$, where $z \in \M$, either
 \begin{enumerate}
   \item Player I wins $G(A)$ with a strategy $\sigma \in \M$; or
   \item Player II wins $G(A)$ with a strategy $\Delta^1_3(z)$-definable over $\M$.
 \end{enumerate}
\end{Theorem}

\begin{Lemma}\label{davislemmalow} Let $z$ be a real and work in $\Pi^1_2(z)\MI + \PPCA$.  Suppose $T \subseteq \om^{<\om}$ is a tree, recursive in $z$, with no terminal nodes, and fix $B \subseteq A \subseteq [T]$ with $B \in \Pi^0_2(z)$ and $A \in \Delta^1_1(z)$.  If $p \in T$ is such that I does not have a winning strategy in $G(A;T_p)$, then there is a quasistrategy $W$ for II in $T_p$ so that
 \begin{itemize}
  \item $[W]\cap B = \emptyset$, and
  \item I does not have a winning strategy in $G(A;W)$.
 \end{itemize}
\end{Lemma}
In keeping with terminology first established in \cite{Da}, we say a position $p$ for which such a quasistrategy $W$ exists is \emph{good} and that $W$ is a \emph{goodness-witnessing quasistrategy} for $p$ (relative to $T,B,A$).

We remark that $\PPCA$ implies $\DCA$, which is equivalent to $\SSAC$ (see VII.6.9 in \cite{Simpson}); this choice principle will be used several times in the course of the proof.

\begin{proof}[Proof of Lemma~\ref{davislemmalow}] Fix a set $U \subseteq \om \times T$ recursive in $z$ so that, setting $U_n = \{ p \in T \mid (n,p) \in U\}$ and $D_n = \{ x \in [T] \mid (\exists k) x \rst k \in U_n\}$, we have $B = \bigcap_{n\in\om} D_n$.  For convenience, we may further assume that each $U_n$ is closed under end-extension in $T$ (i.e., if $p \subseteq q \in T$ and $p \in U_n$, then $q \in U_n$), and that $|p| > n$ whenever $p \in U_n$.

We define an operator $\Phi : \P(T) \to \P(T)$ by setting, for $X \subseteq T$, 
 \begin{align*}
   p \in \Phi(X) \iff (\exists n)&(\forall \sigma) \text{ if }\sigma \text{ is a strategy for I in }T\text{, then}\\
    &(\exists x) x \text{ is compatible with }\sigma, x \notin A\text{, and}
    (\forall k) x\rst k \notin U_n \setminus X.
 \end{align*}
The operator $\Phi$ is clearly monotone on $\P(T)$, and the relation $p \in \Phi(X)$ is $\Pi^1_2(z)$ because this last pointclass is closed under existential quantification over $\omega$ (by $\SSAC$).  We can write this more compactly by introducing an auxiliary game where Player I tries either to get into $A$, or to at some finite stage enter the set $U_n$ while avoiding $X$: Define for $X \subseteq T$ and $n \in \omega$,
 \[
   E^X_n = A \cup \{x \in [T] \mid (\exists k) x\rst k \in U_n \setminus X\}.
 \]
Then
 \[
  p \in \Phi(X) \iff (\exists n) \text{ I doesn't have a winning strategy in }G(E^X_n; T_p).
 \]

Now by $\gP^1_2(z)\MI$ let $\langle \Phi^\alpha \rangle_{\alpha \leq o(\Phi)}$ be the iteration of the operator $\Phi$ with least fixed point $\Phi^\infty$.  (Note this is our sole use of the main strength assumption of the lemma, $\gP^1_2(z)\MI$.)  Let $\phil$ be the associated prewellorder of $\Phi^\infty \subseteq \omega$; formally, we regard definitions and proofs in terms of $\langle \Phi^\alpha \rangle_{\alpha \leq o(\Phi)}$ as being carried out in the theory $\PPCA$ of second order arithmetic, using the real $\phil$ as a parameter.

\begin{Claim} If $p \in T \setminus \Phi^\infty$, then I has a winning strategy in $G(A;T_p)$.
\end{Claim}
\begin{proof}
  For each $q \in T \setminus \Phi^{\infty}$ and $n \in \omega$, we let $\sigma_{q,n}$ be a winning strategy for I in $G(E^{\Phi^\infty}_n;T_q)$, as is guaranteed to exist by the fact that $q \notin \Phi(\Phi^\infty) = \Phi^\infty$.  By $\SSAC$, we may fix a real $\vec{\sigma}$ coding a sequence of such, so that $(\vec{\sigma})_{\la q, n \ra} = \sigma_{q,n}$ for all such pairs $q,n$.

Supposing now that $p \in T \setminus \Phi^\infty$, we describe a strategy $\sigma$ for Player I in $T_p$ from the parameter $\vec{\sigma}$ as follows.  Set $p_0 = p$.  Let $n_0$ be the least $n$ so that $p_0 \notin U_{n}$ (such exists by our simplifying assumption that $|q|>n$ whenever $q \in U_n$).  Suppose inductively that we have reached some position $p_i \notin \Phi^\infty$ and have fixed $n_i$ such that $p_i \notin U_{n_i}$.  Play according to $\sigma_{p_i,n_i}$ until, if ever, we reach a position $q \in U_{n_i} \setminus \Phi^\infty$.  Then set $p_{i+1} = q$, and let $n_{i+1}$ be least such that $p_{i+1} \notin U_{n_{i+1}}$.  

Note the strategy just described is arithmetical in the parameters $z,\vec{\sigma}$, and so exists; call it $\sigma$.  We claim $\sigma$ is winning for I in $G(A;T_p)$.  

Let $x \in [T_p]$ be a play compatible with $\sigma$.  Then $n_0, p_0$ are defined.  If $n_{i+1}$ is undefined for some $i$, then fixing the least such $i$, we must have that no initial segment of $x$ belongs to $U_{n_i} \setminus \Phi^\infty$.  So $x$ is compatible with the strategy $\sigma_{p_i,n_i}$; since this strategy is winning in $G(E^{\Phi^\infty}_{n_i}; T_{p_i})$, we must have that $x \in E^{\Phi^\infty}_{n_i}$.  But then $x \in A$, by definition of the set $E_{n_i}^{\Phi^\infty}$.  

On the other hand, if $n_i$ is defined for all $i$, then by definition of the strategy $\sigma$, we have $p_i \subseteq x$ for all $i$, and for each $i$, $p_i \in \bigcap_{n < n_i} U_n$ (here we use that the sets $U_n$ are closed under end-extension in $T$).  So $x \in \bigcap_{n\in\omega} D_n = B \subseteq A$.  

We have shown $\sigma$ is winning for Player I in $G(A;T_p)$.
\end{proof}
\begin{Claim}
  If $p \in \Phi^\infty$, then $p$ is good.
\end{Claim}
\begin{proof}
  The construction of a quasistrategy $W^p$ witnessing goodness of $p$ proceeds inductively on the ordinal rank of $p \in \Phi^\infty$, that is, on the least $\alpha$ so that $p \in \Phi^\alpha$.  Namely, given such $p, \alpha$, there is some $n$ so that I does not have a winning strategy in the game $G(E^{\Phi^{<\alpha}}_n; T_p)$.  In $W^p$, have II play according to II's non-losing quasistrategy in $G(E^{\Phi^{<\alpha}}_n; T_p)$ until, if ever, a position $q$ in $U_n$ is reached.  Since this non-losing quasistrategy must avoid $U_n\setminus \Phi^{<\alpha}$ by definition of $E_n^{\Phi^{<\alpha}}$, we must have $q \in \Phi^{<\alpha}$; inductively, we have some goodness-witnessing quasistrategy $W^q$ for $q$, so have II switch to play according to this strategy.
  
  Here is a more formal definition of the quasistrategy $W^p$.  For $p \in \Phi^\infty$, define $W^p$ to be the set of positions $q \in T_p$ for which there exists some sequence $\la (\alpha_i, n_i) \ra_{|p| \leq i \leq |q|}$ so that, whenever $|p| \leq i \leq |q|$,
   \begin{itemize}
     \item if $i = |p|$, or $i > |p|$ and $q \rst i \in U_{n_{i-1}}$, then 
       \begin{itemize}
         \item $\alpha_i$ is the least $\alpha$ so that $q \rst i \in \Phi^\alpha$;
         \item $n_i$ is the least $n$ so that I has no winning strategy in $G(E^{\Phi^{<\alpha_i}}_n;T_{q \rst i})$;
       \end{itemize}
     \item if $i > |p|$ and $q \rst i \notin U_{n_{i-1}}$, then $\alpha_i = \alpha_{i-1}, n_i = n_{i-1}$; and
     \item if $i < |q|$, then $q \rst(i+1)$ is in II's non-losing quasistrategy in $G(E^{\Phi^{<\alpha_i}}_{n_i}; T_{q \rst i})$.
   \end{itemize}
Note that formally, we should regard quantification of ordinals $\alpha < o(\Phi)$ as ranging over natural number codes for such as furnished by the prewellorder $\phil$.  The most complicated clauses in the above definition are those involving assertions of the form ``I has no winning strategy in the game $G(E_n^{\Phi^{<\alpha_i}}; T_{q\rst i})$", and such are $\Pi^1_2$ in the parameter $\phil$.  So the criterion for membership in $W^p$ is arithmetical in $\Sigma^1_2(\phil)$ conditions, and therefore by $\PPCA$ the set $W^p$ is guaranteed to exist.  

We need to verify $W^p$ is a quasistrategy for Player II in $T_p$.  An easy induction shows that for each $q \in W^p$, there is a unique witnessing sequence $\la(\alpha_i, n_i)\ra_{|p|\leq i \leq |q|}$ and this sequence depends continuously on $q$; that the $\alpha_i$ are non-increasing; and that I has no winning strategy in $G(E^{\Phi^{<\alpha_i}}_{n_i};T_{q\rst i})$ whenever $|p|\leq i \leq |q|$.  

For $q \in W^p$, we let $\alpha^q, n^q$ denote the final pair (indexed by $|q|$) in the sequence witnessing this membership.  By the above remarks, I has no winning strategy in $G(E^{\Phi^{<\alpha^q}}_{n^q};T_q)$, and by the final condition for membership in $W^p$, the one-step extensions $q^\frown\la l \ra$ in $W^p$ are exactly the one-step extensions of $q$ in II's non-losing quasistrategy in this game.  It follows that $W^p$ is a quasistrategy for II in $T_p$.
%

We claim $W^p$ witnesses goodness of $p$.  We first show $[W^p] \cap B = \emptyset$.  Given any play $x \in [W^p]$, we have some least $i$ so that $\alpha_j = \alpha_i$ for all $j \geq i$; then for all $j>i$, we have $x \rst j$ belongs to II's non-losing quasistrategy in $G(E^{\Phi^{<\alpha_i}}_{n_i}; T_{x \rst i})$.  In particular, for no $k$ do we have $x \rst k \in U_{n_i}$.  Then $x \notin D_{n_i}$, so $x \notin B$ as needed.

We just need to show I has no winning strategy in $G(A; W^p_q)$, for each $q \supseteq p$ in $W^p$.  We argue by induction on $\alpha^q$.  So assume that there is no winning strategy for I in $G(A;W^p_r)$ whenever $\alpha^r < \alpha^q$.

Suppose towards a contradiction that $\sigma$ is a winning strategy for I in $G(A; W^p_q)$.  Let $j$ be least so that $\alpha^q = \alpha_j$.  Then $q$ is in II's non-losing quasistrategy in $G(E^{\Phi^{<\alpha_j}}_{n_j}; T_{q \rst j})$.  We claim no $r \supseteq q$ compatible with $\sigma$ is in $U_{n_j}$.  For otherwise, we have $r \in \Phi^{<\alpha_j}$, so that $\alpha^r < \alpha_j = \alpha^q$, and $\sigma$ is a winning strategy for I in $G(A; W^p_r)$.  This contradicts our inductive hypothesis.

So $\sigma$ cannot reach any position in $U_{n_j}$.  By our definition of $W^p$, we have that the strategy $\sigma$ stays inside II's non-losing quasistrategy for  $G(E^{\Phi^{<\alpha_j}}_{n_j};T_{q\rst j})$.  But since $\sigma$ is winning for I in $G(A;T_{q \rst j})$ and $A \subseteq E^{\Phi^{<\alpha_j}}_{n_j}$, this is a contradiction.

We conclude that I has no winning strategy in $G(A;W^p_q)$; inductively, the claim follows for all $q \in W^p$ extending $p$, so that in particular, $W^p$ witnesses goodness of $p$.
\end{proof}
The last two claims show that every $p \in T$ is either a winning position for I in $G(A;T)$, or is good.  This proves the lemma.
\end{proof}

For future reference, let us refer to the $W^p$ defined in the proof as the \emph{canonical goodness-witnessing strategy for $p$} (relative to $T,B,A$).  We have the following remark, which will be important in computing the complexity of winning strategies:
 \begin{Remark}\label{defcanongood}
   Since ``$\phil$ witnesses the instance of $\Pi^1_2(z)\MI$ at $\Phi$" is $\Delta^1_3(z)$, the statement ``$W$ is the canonical goodness-witnessing strategy for $p$ relative to $T,B,A$" is likewise $\Delta^1_3(z)$ as a relation on pairs $\la W, p \ra$.
 \end{Remark}

\begin{proof}[Proof of Theorem~\ref{detfrommi}]
The proof proceeds from Lemma~\ref{davislemmalow} as usual (see \cite{Da}, \cite{MaBook}); we give a detailed account here, in order to isolate the claimed definability of II's winning strategy.

Fix a $\beta$-model $\M$ of $\PMI$.  Suppose $A$ is $\Sigma^0_3(z)$ for some $z \in \M$; say $A = \bigcup_{k \in \om} B_k$.  By the previous lemma, whenever $T \in \M$ is a tree in $\M$, and $p \in T$ is a position so that in $\M$, there is no winning strategy for I in $G(A; T_p)$, then $p$ is good relative to $T,B_k,A$, for all $k$; that is, for each $k$ there is $W_k$ a quasistrategy for II in $T_p$ so that
  \begin{itemize}
    \item $[W_k] \cap B = \emptyset$;
    \item I does not have a winning strategy in $G(A;W_k)$.
  \end{itemize}  
The idea of the proof is to repeatedly apply the lemma inside $\M$.  At positions $p$ of length $k$, II refines her present working quasistrategy $W_{k-1}$ to one $W_k$ witnessing goodness of $p$ relative to $W_{k-1},B_k,A$, so ``dodging" each of the $\Pi^0_2(z)$ sets $B_k$, one at a time. 

More precisely: suppose I does not win $G(A)$ in $\M$, where $A$ is $\Sigma^0_3(z)$ for some $z \in \M$.  Let $W^\emptyset$ be the canonical goodness-witnessing quasistrategy for $\emptyset$ relative to $\om^{<\om}, B_0, A$ as constructed in the proof of Lemma~\ref{davislemmalow}.  Then let $H^\emptyset$ be II's non-losing quasistrategy in $G(A;W^\emptyset)$ (so that for no $p \in H^\emptyset$ do we have that I wins $G(A;W^\emptyset_p)$).

Suppose inductively that for some $k$, we have subtrees $H^p$ of $T$, defined for a subset of $p \in T$ with length $\leq k$, so that
 \begin{enumerate}
   \item each $H^p$ is a quasistrategy for II in $T_p$ and belongs to $\M$;
   \item $[H^p] \cap B_{|p|} = \emptyset$;
   \item for no $q \in H^p$ does I have a winning strategy in $G(A; H^p_q)$;
   \item if $p \subseteq q$, then $H^q \subseteq H^p$ whenever both are defined;
   \item if $|p|<k$ and $p^\frown \la l \ra \in H^p$, then $H^{p^\frown \la l \ra}$ is defined.
 \end{enumerate}
In order to continue the construction, we need to define quasistrategies $H^{p^\frown \la l \ra}$, whenever $|p|=k$, $H^p$ is defined, and $p^\frown{\la l \ra} \in H^p$.  Given such $p$ and $l$, we have that I has no winning strategy in $G(A;H^p_{p^\frown\la l \ra})$ by (3).  So applying Lemma~\ref{davislemmalow} inside $\M$, let $W^{p^\frown \la l \ra}$ be the canonical goodness-witnessing strategy for $p^\frown \la l \ra$ relative to $H^p, B_{k+1}, A$.  Then let $H^{p^\frown\la l \ra}$ be II's non-losing quasistrategy in $G(A;W^{p^\frown\la l \ra})$.  It is easy to see that this quasistrategy satisfies the properties (1)-(4), so we have the desired system of quasistrategies $H^q$ satisfying (5), for $|q|=k+1$.

Now set $p \in H$ if and only if for all $i < |p|$, $H^{p \rst i}$ is defined and $p \in H^{p \rst i}$.  It follows from (5) that $H$ is a quasistrategy for II, and by (4) we have $H \subseteq H^p$ for each $p \in H$.  By (2) then, $[H]\cap B_k = \emptyset$ for all $k \in \om$, so that $[H] \cap A = \emptyset$.

Observe that for each $p\in H$, we have that the sequence $\la H^{p \rst i} \ra_{i < |p|}$ exists in $\M$, since it is obtained by a finite number of applications of $\PMI$ and $\PPCA$.  Since $\M$ is a $\beta$-model, it really is the case (in $V$) that $[H] \cap B_k = \emptyset$ for all $k \in \om$.  Though $H$ need not belong to $\M$, we claim it is nonetheless a $\Delta^1_3(z)$-definable class over $\M$.  For $p \in H$ if and only if there exists a sequence $\la W_i, H_i\ra_{i < |p|}$, so that for all $i < |p|$,
 \begin{itemize}
   \item $W_i$ is II's canonical goodness-witnessing strategy for $p \rst i$, relative to $H_{i-1},B_i,A$ (where we set $H_{-1}=\om^{<\om}$);
   \item $H_{i}$ is II's non-losing quasistrategy in $G(A;W_i)$ at $p \rst i$;
   \item for all $i < |p|$, $p \in H_i$.
 \end{itemize}
This is a $\Sigma^1_3(z)$ condition, by Remark~\ref{defcanongood}.  And note that $p \notin H$ if and only if there is a sequence $\la H_i, W_i\ra_{i \leq l}$, for some $l < |p|$, satisfying the first two conditions for $i\leq l$, but so that $p \notin H_{l}$.  This is likewise $\Sigma^1_3(z)$, so that $H$ is $\Delta^1_3(z)$-definable in $\M$.

Given a $\Delta^1_3(z)$ definition of the quasistrategy $H$, it is easy to see that the strategy $\tau$ for II obtained by taking $\tau(p)$ to be the least $l$ so that $p^\frown{\la l \ra} \in H$ is likewise $\Delta^1_3(z)$ and winning for II in $G(A;T)$.  This completes the proof of Theorem~\ref{detfrommi}.
\end{proof}

\section{$\boldsymbol{\Pi^1_2}$ monotone induction from infinite depth $\Sigma_2$-nestings}\label{sec:pmi}
In this section, the theories of $\KP$ and $\Sigma_1$-Comprehension are defined in the language of set theory as usual.  We will furthermore make use of the theories $\KPI_0$, which asserts that every set is contained in some admissible set (that is, some transitive model of $\KP$), and $\KPI$, which is the union of $\KP$ and $\KPI_0$.  $\KPI_0$ is relevant largely because it is a weak theory in which Shoenfield absoluteness holds; in particular, $\boldsymbol{\Pi^1_2}$ expressions are equivalent over $\KPI_0$ to $\Pi_1$ statements in the language of set theory.

We remark that $\KPI_0$ and $\bgP^1_1\CA_0$ prove the same statements of second order arithmetic.  Since we primarily work with models in the language of set theory in this section, we take $\KPI_0$ as our base theory, but all of the results proved here can be appropriately reformulated as statements about countably coded $\beta$-models in second order arithmetic (as in Chapter VII of \cite{Simpson}).

For $\M$ an illfounded model in the language of set theory, we identify the wellfounded part of $\M$ with its transitive collapse, denote this $\wfp(\M)$, and set $\wfo(\M) = \wfp(\M) \cap \ON$.  Recall we say $\M$ is an $\omega$-model if $\omega < \wfo(\M)$.  The following definition is due to Welch \cite{We2}.

\begin{Definition}\label{IDN} For $\M$ an illfounded $\omega$-model of $\KP$ in the language of set theory, an \emph{infinite depth $\Sigma_2$-nesting based on $\M$} is a sequence $\langle \zeta_n, s_n \rangle_{n \in \om}$ of pairs so that for all $n \in \om$,
 \begin{enumerate}
  \item $\zeta_n \leq \zeta_{n+1} < \wfo(\M)$,
  \item $s_n \in \ON^{\M} \setminus \wfo(\M)$,
  \item $\M \models s_{n+1} < s_n$,
  \item $(L_{\zeta_n} \prec_{\Sigma_2} L_{s_n})^{\M}$.
 \end{enumerate}
\end{Definition}

\begin{Lemma} Suppose $\gamma_1 \leq \gamma_2 < \delta_2 < \delta_1$ are ordinals so that
 \begin{enumerate}
   \item $L_{\gamma_1} \prec_{\Sigma_1} L_{\delta_1}$;
   \item $L_{\gamma_2} \prec_{\Sigma_2} L_{\delta_2}$;
   \item $\delta_1$ is the least admissible ordinal above $\delta_2$;
   \item For all $\alpha \leq \delta_2$, $L_\alpha$ $\Sigma_\omega$-projects to $\omega$.
 \end{enumerate}
Then $L_{\gamma_2}$ satisfies $\Pi^1_2(z)\MI$, for all reals $z \in L_{\gamma_1}$.
\end{Lemma}
Item (4) simply asserts that for every $\alpha$, there is a subset of $\omega$ definable over $L_\alpha$ that doesn't belong to $L_\alpha$; this simplifying assumption ensures that every $L_\alpha$ is countable, as witnessed by a surjection $f:\om\to L_\alpha$ that belongs to $L_{\alpha+1}$.  Note the least level of $L$ that does \emph{not} $\Sigma_\omega$-project to $\omega$ is a model of $\ZF^{-}$; since this is far beyond the strength of the theories considered here, we don't lose anything by assuming (4).
\begin{proof}
Let $\Phi: \P(\om) \to \P(\om)$ be a $\Pi^1_2(z)$ monotone operator in $L_{\gamma_2}$; fix a $\Pi^0_1(z)$ condition $T$ so that
 \[
  n \in \Phi(X) \iff L_{\gamma_2} \models (\forall x)(\exists y)T(n,X,x,y,z)
 \]
whenever $n \in \om$ and $X \in \P(\om) \cap L_{\gamma_2}$.  Notice that for such $X$
 \[
  n \in \Phi(X) \iff (\forall x \in L_{\gamma_2})(\exists y)T(n,X,x,y,z),
 \]
by absoluteness and because $\gamma_2$ is a limit of admissibles.  Regarding the operator defined in this way, ``$n \in \Phi(X)$" makes sense even for sets $X \notin L_{\gamma_2}$ (though this extended $\Phi$ may fail to be monotone when applied to sets not in $L_{\gamma_2}$).

For each ordinal $\eta$, we define the approximation $\Phi_\eta$ as the operator $\Phi$ relativized to $L_\eta$,
 \[
  n \in \Phi_\eta(X) \iff (\forall x \in L_\eta)(\exists y)T(n,X,x,y,z).
 \]
The point is that the operator $\Phi_\eta$ is then $\Sigma^1_1$ in any real parameter coding the countable set $\R \cap L_\eta$ (for example, $\Th(L_\eta)$, the characteristic function of the theory of $L_\eta$ under some standard coding), and so each $\Phi_\eta$ will be correctly computed in, e.g., $L_\alpha$ for $\alpha$ a limit of admissibles above $\eta$.

Obviously $\Phi = \Phi_{\gamma_2}$, so is monotone in $L_{\gamma_2}$.  But for $\eta \neq \gamma_2$ we may not even have that the operators $\Phi_\eta$ are monotone on $\P(\om) \cap L_\eta$.  So we instead work with the obvious ``monotonizations",
 \begin{align*}
  n \in \Psi_\eta(X) &\iff (\exists X' \subseteq X) n \in \Phi_{\eta}(X') \\ &\iff (\exists X' \subseteq X)(\forall x \in L_{\eta})(\exists y)T(n,X',x,y,z).
 \end{align*}
These are again $\Sigma^1_1(\Th(L_\eta), z)$, and $\Psi_{\gamma_2}(X) = \Phi_{\gamma_2}(X) = \Phi(X)$ for $X \in L_{\gamma_2}$.

Let $\la \Psi^\xi_\eta \ra_{\xi \leq o(\Psi_\eta)}$ be the sequence obtained via iterated application of the operator $\Psi_\eta$, as in Definition~\ref{def:indops}.  The most important properties of these sequences are captured in the following two claims.

\begin{Claim}
 If $\eta < \eta'$, then $(\forall X)\Psi_\eta(X) \supseteq \Psi_{\eta'}(X)$.
\end{Claim}
\begin{proof}
  Suppose $n \in \Psi_{\eta'}(X)$; then
   \[
    (\exists X' \subseteq X)(\forall x \in L_{\eta'})(\exists y)T(n,X',x,y,z),
   \]
  and any such $X'$ will likewise be a witness to $n \in \Psi_{\eta}(X)$, since the latter is defined the same way but with the universal quantifier bounded by the smaller set $L_\eta$.
\end{proof}

\begin{Claim} Suppose $\xi < \xi'$ and $\eta < \eta'$.  Then
  \begin{enumerate}
   \item $\Psi_\eta^\xi \subseteq \Psi_\eta^{\xi'}$;
   \item $\Psi_\eta^\xi \supseteq \Psi_{\eta'}^{\xi}$.
  \end{enumerate}
\end{Claim}
\begin{proof}
  (1) is by definition.  (2) follows from induction and the chain of inclusions, for $X \supseteq Y$,
   \[
     \Psi_{\eta'}(Y) \subseteq \Psi_{\eta'}(X) \subseteq \Psi_\eta(X),
   \]
 the first by monotonicity of $\Psi_{\eta'}$, the second by the previous claim.
\end{proof} 
So the array $\la \Psi_\eta^\xi \ra$ is increasing in $\xi$ and decreasing in $\eta$.  Applying this claim with $\xi = \omega_1$, we have $\Psi_\eta^\xi = \Psi_\eta^\infty$, so that $\Psi_\eta^\infty \supseteq \Psi_{\eta'}^\infty$ whenever $\eta < \eta'$.  

We now consider definability issues with respect to the operators $\Psi_\eta$ and the associated sequences, with the aim of showing the levels of $L$ under consideration are sufficiently closed to correctly compute these objects, and ultimately ensuring that the sequences $\la \Psi^\xi_\eta \ra_{\xi \leq o(\Psi_\eta)}$ converge to the sequence of interest $\la \Psi^\xi_{\gamma_2} \ra_{\xi \leq o(\Psi_{\gamma_2})}$ as $\eta \to \gamma_2$.

Note that the hypothesis $L_{\gamma_2} \prec_{\Sigma_2} L_{\delta_2}$ implies $L_{\gamma_2}$ is a model of $\Sigma_2$-$\KP$.  The assumed elementarity in the lemma then implies each of $\gamma_1,\gamma_2,\delta_2$ is a limit of $\Sigma_2$-admissible ordinals.

\begin{Claim}\label{seqsinside} Suppose $z \in L_\alpha$ and $L_\alpha \models \KPI$.  Then the relation ``$n \in \Psi^\xi_\eta$" (as a relation on $\la n, \xi, \eta \ra \in \om \times \alpha \times \alpha$) is $\Delta_1^{L_\alpha}$ in the parameter $z$.  Consequently, for all $\eta < \alpha$ and $\nu < \alpha$, the sequence $\la \Psi^\xi_\eta \ra_{\xi < \nu}$ belongs to $L_\alpha$.
\end{Claim}
\begin{proof}
 The relation $n \in \Psi_\eta(X)$ is, as remarked above, $\Sigma^1_1(z, \Th(L_\eta))$ on $n,X$, and so is $\Pi_1$ over the least admissible set containing $z, \eta$.  Since every set is contained in some admissible set $L_\beta$ with $\beta<\alpha$, we have that ``$n \in \Psi_\eta(X)$" is $\Delta_1(z)$ over $L_\alpha$.  The last part of the claim then follows from $\Sigma_1$-recursion inside $L_\alpha$, using the $\Delta_1^{L_\alpha}(z)$-definability of the relation $Y = \Psi_\eta(X)$.
 \end{proof}

\begin{Claim} Suppose $z \in L_\alpha$, $\eta < \alpha$ and $L_\alpha$ is a model of $\Sigma_1$-Comprehension.  Then $o(\Psi_\eta) < \alpha$, and $\la \Psi_\eta^\xi \ra_{\xi \leq o(\Psi_\eta)} \in L_\alpha$.  Moreover, the relation $n \in \Psi_\eta^\infty$ (on $\omega \times \alpha$) is $\Delta_1^{L_\alpha}(z)$.
\end{Claim}
\begin{proof}  Note such $L_\alpha$ satisfies $\KPI$, so by the previous claim together with $\Sigma_1$-Comprehension, $P_\eta = \{n\in\om \mid (\exists \xi<\alpha) n \in \Psi_\eta^\xi\} \in L_\alpha$.  By admissibility, the map on $P_\eta$ sending $n$ to the least $\xi$ such that $n \in \Psi_\eta^\xi$ is bounded in $\alpha$, and the claim is immediate.  The last assertion holds because in $L_\alpha$,
 \[
  n \in \Psi_\eta^\infty \iff (\exists \xi) n \in \Psi^\xi_\eta \iff (\forall \xi)(\Psi_\eta^\xi = \Psi_\eta^{\xi+1} \to n \in \Psi_\eta^\xi).
 \]
\end{proof}
\begin{Claim} Suppose $z \in L_\alpha$, and that $\alpha$ is a limit of ordinals $\beta$ so that $L_\beta$ is a model of $\Sigma_1$-Comprehension.  Then the relation $n \in \Psi_\eta^\infty$ is $\Delta_1^{L_\alpha}(z)$.
\end{Claim}
\begin{proof} Immediate from the previous claim and the fact that the sequences are correctly computed in models of $\KPI_0$.
\end{proof}

\begin{Claim}\label{stabcolumns} If $\xi < \gamma_2$, then for some $\eta_0 < \gamma_2$ we have $\Psi_{\eta_0}^\xi = \Psi_{\gamma_2}^\xi$; furthermore, $\la \Psi_{\gamma_2}^\zeta \ra_{\zeta < \xi} \in L_{\gamma_2}$.
\end{Claim}
Note that then for this $\eta_0$, $\Psi_{\eta_0}^\xi = \Psi_\eta^\xi$ whenever $\eta_0 \leq \eta < \delta_2$.
\begin{proof}
The set $Q_\xi=\{n\in\om \mid (\exists \eta<\gamma_2) n \notin \Psi_\eta^\xi\}$ is a member of $L_{\gamma_2}$ by $\Sigma_1$-Comprehension there.  Now the map sending $n \in Q_\xi$ to the least $\eta$ such that $n \notin \Psi^\xi_\eta$ is $\Delta_1$, so by admissibility, is bounded by some $\eta_0 < \alpha$.  Recall the sequence $\la \Psi^\xi_\eta\ra_{\eta \in \ON}$ is decreasing in $\eta$; so 
 \begin{align*}
   n \in \Psi^\xi_{\eta_0} \iff L_{\gamma_2} \models (\forall \eta) n \in \Psi^\xi_\eta
   &\iff L_{\delta_2} \models (\forall \eta) n \in \Psi_\eta^\xi \\
   &\Longrightarrow n \in \Psi^\xi_{\gamma_2} \Longrightarrow n \in \Psi^\xi_{\eta_0}.
 \end{align*}
Note we have used the fact that $L_{\gamma_2} \prec_{\Sigma_1} L_{\delta_2}$.  For the last part of the claim, consider the map sending $\zeta < \xi$ to the least $\eta_0$ such that $(\forall \eta > \eta_0) \Psi_\eta^\zeta = \Psi_{\eta_0}^\zeta$.  This map is $\Pi_1$-definable, so by $\Sigma_2$-Collection in $L_{\gamma_2}$, we have a bound $\bar{\eta} < \gamma_2$, and for each $\zeta < \xi$, $\Psi^\zeta_{\bar{\eta}} = \Psi^\zeta_{\gamma_2}$.  By Claim~\ref{seqsinside} the sequence $\la \Psi_{\bar{\eta}}^\zeta \ra_{\zeta < \xi} = \la \Psi_{\gamma_2}^\zeta \ra_{\zeta < \xi}$ is in $L_{\gamma_2}$.
\end{proof}

\begin{Claim} For all $\xi < \gamma_2$, $\Psi^\xi_{\gamma_2} = \Psi^\xi_{\delta_2}$; consequently, $\Psi_{\gamma_2}^{<\gamma_2} = \Psi_{\delta_2}^{<\gamma_2}$.
\end{Claim}
\begin{proof}
By using induction on $\xi$ and since $\Psi^{<\xi}_{\gamma_2} \in L_{\gamma_2}$ by the previous claim, it is sufficient to show $\Psi_{\gamma_2}(X) = \Psi_{\delta_2}(X)$ whenever $X \in L_{\gamma_2}$.  We already know $\supseteq$ holds.

So suppose $n \in \Psi_{\gamma_2}(X)$.  Then we have $n \in \Phi_{\gamma_2}(X) = \Phi(X)$, by monotonicity of $\Phi = \Phi_{\gamma_2}$ in $L_{\gamma_2}$.  So
  \[ 
   L_{\gamma_2} \models (\forall x)(\exists y) T(n,X,x,y,z)
  \]
so that by $\Sigma_1$-elementarity (this is enough, since $\Pi^1_2$ relations are $\Pi^{\KPI_0}_1$), $L_{\delta_2}$ models the same.  Thus $n \in \Psi_{\delta_2}(X)$ (with witness $X' = X$).
\end{proof}
We haven't yet used the full strength of $L_{\gamma_2} \prec_{\Sigma_2} L_{\delta_2}$, nor, for that matter, any of the assumptions on $\gamma_1, \delta_1$.  We appeal to the first assumption to show that in fact $o(\Psi_{\delta_2}) \leq \gamma_2$; the second will be used to show that $\Psi_{\delta_2}^\infty = \Psi_{\gamma_2}^\infty$, and it will follow that the operator $\Psi_{\gamma_2}$ (which is equal to $\Phi$, remember) stabilizes inside $L_{\gamma_2}$.

Notice that by Claim~\ref{stabcolumns}, $\Psi_{\gamma_2}^\xi = \bigcap_{\eta < \gamma_2} \Psi_\eta^\xi$ for all $\xi < \gamma_2$.  So
 \[
   \Psi_{\gamma_2}^{<\gamma_2} = \{ n \in \om \mid (\exists \xi < \gamma_2)(\forall \eta < \gamma_2) n \in \Psi_\eta^\xi\}.
 \]
This set is $\Sigma_2$-definable over $L_{\gamma_2}$.  By the fact that $\Psi_{\delta_2}^{\gamma_2} \subseteq \Psi_{\eta}^{\gamma_2}$ for all $\eta < \delta_2$, we have
 \[
   \Psi_{\delta_2}^{\gamma_2} \subseteq \{n \in \om \mid (\forall \eta< \delta_2) n \in \Psi^{\gamma_2}_\eta \} \subseteq \{n \in \om \mid (\exists \xi < \delta_2)(\forall \eta < \delta_2) n \in \Psi^\xi_\eta\}.
 \]
By the assumed $\Sigma_2$-elementarity $L_{\gamma_2} \prec_{\Sigma_2} L_{\delta_2}$, this last set is precisely $\Psi_{\gamma_2}^{<\gamma_2}$.  We obtain
 \[
  \Psi_{\delta_2}^{\gamma_2} \subseteq \Psi_{\gamma_2}^{<\gamma_2} = \Psi_{\delta_2}^{<\gamma_2} \subseteq \Psi_{\delta_2}^{\gamma_2}
 \]
so that $\Psi^{\gamma_2}_{\delta_2} = \Psi^{<\gamma_2}_{\delta_2}$ is the least fixed point of $\Psi_{\delta_2}$, $\Psi_{\gamma_2}^{<\gamma_2} = \Psi_{\delta_2}^\infty$.

\begin{Claim} $\Psi_{\delta_2}^\infty = \Psi_{\gamma_2}^\infty$.
\end{Claim}
\begin{proof}  
As usual, we know $\subseteq$ holds since $\gamma_2 < \delta_2$.  We have $\Psi_{\delta_2}^\infty = \Psi_{\gamma_2}^{<\gamma_2} \in L_{\delta_2}$.  Suppose $n \notin \Psi_{\delta_2}^\infty$.  Then
 \[
  L_{\delta_1} \models (\exists \eta)(\exists P)(\forall m \in \omega)(m \in \Psi_\eta(P) \to m \in P) \wedge n \notin P,
 \]
with $\eta = \delta_2$ and $P = \Psi_{\delta_2}^\infty$. Recall ``$m \in \Psi_\eta(P)$", being a $\Sigma^1_1$ statement about $m, \Th(L_\eta), P$, is $\Pi_1$ over any admissible set containing $\eta,z,P$.  Since $L_{\delta_1}$ is assumed to be admissible, the relation above is then $\Sigma_1$ in $L_{\delta_1}$.  It therefore reflects to $L_{\gamma_1}$ (recall that $z$, the parameter from which everything is defined, is assumed to belong to $L_{\gamma_1}$).  But then $n \notin \Psi_\eta^\infty$ for some $\eta < \gamma_1$; hence $n \notin \Psi_{\gamma_2}^\infty$.
\end{proof}

So the least fixed points $\Phi^\infty = \Psi_{\gamma_2}^\infty$ and $\Psi_{\delta_2}^\infty$ are equal.  The argument just given shows the relation $n \notin \Phi^\infty$ is $\Sigma_1$ over $L_{\delta_1}$, hence over $L_{\gamma_1}$; in any event, the set $\Phi^\infty$ belongs to $L_{\gamma_2}$ (using $\Sigma_1$-Comprehension in $L_{\gamma_2}$ in the case that $\gamma_1=\gamma_2$).

Finally, we claim $o(\Phi) < \gamma_2$.  The map defined in $L_{\gamma_2}$ that takes $n \in \Phi^\infty = \Psi_{\gamma_2}^\infty$ to the least $\xi$ such that $(\exists \eta_0)(\forall \eta > \eta_0) n \in \Psi_{\eta}^\xi$ is $\Sigma_2$-definable, and so by $\Sigma_2$-Collection is bounded in $\gamma_2$.  Since for each $\xi < \gamma_2$ we have $\Phi^\xi = \Psi_{\gamma_2}^{\xi} = \Psi_{\eta_0}^\xi$ for some $\eta_0 < \gamma_2$, this implies $o(\Phi) < \gamma_2$.

That $\la \Phi^\xi \ra_{\xi \leq o(\Phi)}$ belongs to $L_{\gamma_2}$ now follows from the last assertion of Claim~\ref{stabcolumns}.  This completes the proof that the desired instance of $\Pi^1_2(z)\MI$ holds in $L_{\gamma_2}$.
\end{proof}

\begin{Theorem}\label{PMIfromIDN}
Suppose $\M$ is an illfounded $\omega$-model of $\KP$ with $\la \zeta_n, s_n \ra_{n \in \omega}$ an infinite depth $\Sigma_2$-nesting based on $\M$, and that $\M$ is locally countable, in the sense that every $L_a^{\M}$ has ultimate projectum $\omega$ in $\M$.  Then if $\beta = \sup_{n\in\om} \zeta_n$, we have $L_\beta \models \PMI$.
\end{Theorem}
\begin{proof} If $\beta= \zeta_n$ for some $n \in \omega$, then we obtain the result immediately by applying the lemma in $M$ to the tuple $\langle \zeta_n, \zeta_{n+1}, s_{n+1},s_n \rangle$.  So we can assume $\la \zeta_n \ra_{n \in \om}$ is strictly increasing.  Let $\Phi : \P(\om) \to \P(\om)$ be $\Pi^1_2(z)$ and monotone in $L_\beta$ for some $z \in L_\beta$, and let $\zeta_n$ be sufficiently large that $z \in L_{\zeta_n}$.  Now $L_{\zeta_{n+1}} \prec_{\Sigma_1} L_\beta$ and both models satisfy $\KPI_0$, so that whenever $X \subseteq \om$ is in $L_{\zeta_{n+1}}$, we have
 \[
  L_\beta \models n \in \Phi(X) \iff L_{\zeta_{n+1}} \models n \in \Phi(X).
 \]
In particular, $L_{\zeta_{n+1}}$ believes $\Phi$ is $\Pi^1_2(z)$ and monotone, so that by the lemma applied to the tuple $\la \zeta_n, \zeta_{n+1}, s_{n+1}, s_n\ra$, we have $o(\Phi) < \zeta_{n+1}$, and the sequence $\la \Phi^\xi \ra_{\xi \leq o(\Phi)}$ (which is computed identically in $L_{\zeta_{n+1}}$ and $L_\beta$) belongs to $L_{\zeta_{n+1}}$.
\end{proof}
Combining Theorems~\ref{detfrommi} and \ref{PMIfromIDN}, we obtain
\begin{Corollary} If $\M,\beta$ are as in the previous theorem, then for any $\gS^0_3(z)$ set with $z \in L_\beta$, either
  \begin{enumerate}
   \item Player I wins $G(A)$ with a strategy $\sigma \in L_\beta$; or
   \item Player II wins $G(A)$ with a strategy $\Delta^1_3(z)$-definable over $L_\beta$.
 \end{enumerate}
\end{Corollary}

\section{Infinite depth $\Sigma_2$-nestings from determinacy}\label{sec:wel}
In this section we show in the base theory $\PCA$ that $\Sigma^0_3$-DET implies the existence of models bearing infinite depth $\Sigma_2$ nestings.  The arguments are mostly cosmetic modifications of those given in Welch's \cite{We}.  The most significant adjustment is to the Friedman-style game, Welch's $G_\psi$, which is here tailored to allow the proof of the implication to be carried out in $\PCA$.

For $\alpha$ an ordinal, let $T^\alpha_2$ denote the lightface $\Sigma_2$-theory of $L_\alpha$, i.e.,
 \[
  T^\alpha_2 = \{ \sigma \mid \sigma \text{ is a }\Sigma_2\text{ sentence without parameters, and }L_\alpha \models \sigma\}.
 \]
We will also abuse this notation slightly by applying it to nonstandard ordinals $b$, so that if $b \in \ON^{\M} \setminus \wfo(\M)$, $T^b_2$ denotes the $\Sigma_2$-theory of $(L_b)^{\M}$.  It will always be clear from context which illfounded model $\M$ this $b$ comes from.

\begin{Lemma}\label{welchlem}
  Suppose $\M$ is an illfounded $\omega$-model of $\KP$ such that $(L_a)^{\M} \models ``$all sets are countable", for every $a \in \ON^{\M}$.  Set $\beta = \wfo(\M)$.  Suppose for all nonstandard ordinals $a$ of $\M$, there exists some $<^{\M}$-smaller nonstandard $\M$-ordinal $b$ so that $T^b_2 \subseteq T^\beta_2$.  Then there is an infinite depth $\Sigma_2$ nesting based on $\M$.
\end{Lemma}

\begin{proof}
  This is essentially shown in Claim (5) in Section 3 of \cite{We}.  We outline the shorter approach suggested there.
  
  Suppose $b$ is a nonstandard $\M$-ordinal with $T^b_2 \subseteq T^\beta_2$.  By the assumption of local countability in levels of $L^\M$, we have a uniformly $\Sigma_2$-definable $\Sigma_2$ Skolem function, which we denote $h^b_2$ (see \cite{SDFrUnif}).  The set $H = h_2^b [\om^{<\om}]$ is transitive in $\M$, since for any $x \in H$, the $<_L^{\M}$-least surjection of $\omega$ onto $x$ is in $H$, and since $\M$ is an $\omega$-model, the range of this surjection is a subset of $H$.  Since $H \models V=L$, we have by condensation in $\M$ that $H = L_{\gamma_b} \prec_{\Sigma_2} L_b$ for some $\gamma_b \leq^{\M} b$.
  
  We claim that $\gamma_b < \beta$.  For suppose not, so there is some nonstandard ordinal $c$ of $L_b$ in $L_{\gamma_b}$.  Let $f$ be the $<_L^{\M}$-least surjection from $\om$ onto $c$.  Then $f = h^b_2(k)$ for some $k \in \om$, and for $m,n\in \om$, the sentences ``$h_2(k)$ exists, is a function from $\omega$ onto some ordinal, and $h_2(k)(m) \in h_2(k)(n)$" are $\Sigma_2$.  But since $T^b_2 \subseteq T_2^\beta$, this would imply $h^\beta_2(k)(m) \in h_2^\beta(k)(n)$ whenever $f(m) \in f(n)$ in $(L_b)^\M$.  This contradicts the wellfoundedness of $\beta$.
  
  The lemma now follows by choosing some descending sequence $\la b_n \ra_{n \in \om}$ of nonstandard ordinals of $\M$ with $T^{b_n}_2 \subseteq T^{\beta}_2$ for all $n$, and setting $\gamma_n = \sup h_2^{b_n}[\om^{<\om}] < \beta$.  Since the $\gamma_n$ are true ordinals, we can choose some non-decreasing subsequence $\la \gamma_{n_k} \ra_{k \in \om}$, and $\la \gamma_{n_k}, b_{n_k} \ra_{k \in \om}$ is the desired infinite depth $\Sigma_2$-nesting.
\end{proof}

\begin{Theorem}\label{IDNfromDet} Work in $\PCA$.  If $\Sigma^0_3$-determinacy holds, then there is a model $\M$ for which there exists an infinite depth $\Sigma_2$-nesting based on $\M$.
\end{Theorem}
\begin{Corollary} Work in $\PCA$.  $\Sigma^0_3$-determinacy implies the existence of a $\beta$-model of $\PMI$; indeed, $L_\gamma \models \PMI$ for some countable ordinal $\gamma$.  
\end{Corollary}
\begin{proof} Immediate, combining Theorem~\ref{IDNfromDet} with Theorem~\ref{PMIfromIDN}.  \end{proof}
\begin{proof}[Proof of Theorem~\ref{IDNfromDet}]  We define a variant of Welch's game $G_\psi$ from \cite{We}.  Players I and II play complete consistent theories in the language of set theory, $\fone, \ftwo$, respectively, extending 
 \begin{equation*}\tag{$*$}\label{VLprojom}
   V=L\;\; +\;\; KP\;\; + \;\;\rho_\om=\om.
 \end{equation*}
These theories uniquely determine term models $\Mi, \Mii$.  Player I loses if $\Mi$ has nonstandard $\omega$; similarly, if $\Mi$ is an $\om$-model and $\Mii$ is not, then Player II loses.  (Note that this is a Boolean combination of $\Sigma^0_2$ conditions on $\fone, \ftwo$.)

The remainder of the winning condition assumes $\Mi, \Mii$ are both $\omega$-models.  Player I wins if any of the following hold.
 \begin{enumerate}
   \item $\ftwo \in \Mi$, or $\fone = \ftwo$.
   \item $(\exists \beta \leq \ON^{\Mi})(\exists a \in \ON^{\Mii})(\forall n \in \om)(\exists \la a_i, s_i \ra_{i \leq n})$ so that, for all $i < n$,
     \begin{itemize}
      \item $a_0 = a$ and $a_i \in \ON^{\Mii}$,
      \item $(a_{i+1} < a_i)^{\Mii}$,
      \item $\sigma_i$ is the first $\Sigma_2$ formula (in some fixed recursive list of all formulas in the language of set theory) so that $L_\beta^{\Mi} \not\models \sigma_i$ and $L_{a_i}^{\Mii} \models \sigma_i$;
      \item if $a_i$ is a successor ordinal in $\Mii$, then $a_{i+1}$ is the largest limit ordinal of $\Mii$ below $a_i$;
      \item if $a_i$ is a limit ordinal in $\Mii$, and $\sigma_i$ is the formula $\exists u \forall v \psi(u,v)$, then $a_{i+1}$ is least in $\ON^{\Mii}$ so that $(\exists u \in L_{a_{i+1}})(L_{a_i} \models \forall v \psi(u,v))$ in $\Mii$.
     \end{itemize}
 \end{enumerate}
Note that if (2) holds, then $\Mii$ must be an illfounded model: If $\beta, a$ witness the condition, then the sequences $\la a_i, \sigma_i \ra_{i \leq n}$ are uniquely determined for each $n$, and are inclusionwise increasing in $n$; then $\la a_i \ra_{i \in \om}$ is an infinite descending sequence of $\Mii$-ordinals.

Note also that this condition is $\Sigma^0_3$ as a condition on $\fone, \ftwo$.  Strictly speaking, the quantifiers over $\Mi$, $\ON^{\Mii}$, etc. should be regarded as natural number quantifiers ranging over the indices of defining formulas for members of the models $\Mi,\Mii$.  Clause (1) is then $\Sigma^0_2$, and (2) is $\Sigma^0_3$, since each bulleted item there is recursive in codes for the objects $\beta,a,\la a_i, \sigma_i \ra_{i \leq n}$ and the pair $\la \fone, \ftwo \ra$.

Denote the set of runs which I wins by $F$; so $F$ is $\Sigma^0_3$.

\begin{Claim} Player I has no winning strategy in $G(F)$.
\end{Claim}
\begin{proof} Suppose instead that I has some winning strategy in this game.  By Shoenfield absoluteness (which holds in $\PCA$, see \cite{Simpson}) there is such a winning strategy $\sigma$ in $L$.  Let $\alpha$ be the least admissible ordinal so that $\sigma \in L_\alpha$ (such exists since $\PCA$ implies the reals are closed under the hyperjump; see \cite{Sacks}).  Let $\ftwo$ be the theory of $L_\alpha$.  Note that then $L_\alpha$ $\Sigma_1$-projects to $\omega$, since it is the least admissible containing some real; in particular, it satisfies condition (\ref{VLprojom}).  Let $\fone = \sigma * \ftwo$ be the theory that $\sigma$ responds to $\ftwo$ with.

Now $\sigma$ is winning for I in $G(F)$; so $\Mi$ is an $\om$-model.  Since $\Mii$ is wellfounded, (2) must fail, and since we assumed $\sigma$ is winning for I, we have (1) holds; that is, either $\ftwo \in \Mi$ or $\fone = \ftwo$.  If $\fone = \ftwo$, then II was simply copying I's play, so that $\sigma \in L_\alpha = \Mi$, implying $\fone \in \Mi$, a contradiction to the fact that $\Mi$ $\omega$-projects to $\om$.

So $\ftwo \in \Mi$.  The strategy $\sigma$ is computable from $\ftwo$, so must also belong to $\Mi$.  But then, since $\fone = \sigma * \ftwo$, we again obtain the contradiction $\fone \in \Mi$.
\end{proof}

\begin{Claim} If there is no model with an infinite depth $\Sigma_2$-nesting, then Player II has no winning strategy in $G(F)$.
\end{Claim}
\begin{proof} Towards a contradiction, let $\tau$ be a winning strategy for II; as in the previous claim, we may assume $\tau \in L$, and let $\alpha$ be the least admissible with $\tau \in L_\alpha$.  Put $\fone = \Th(L_\alpha)$; then $\fone$ satisfies the condition (\ref{VLprojom}).  Let $\ftwo = \tau * \fone$ be $\tau$'s response.

We claim that if $\Mii$ is the model so obtained, then $\wfo(\Mii) \leq \alpha$ (note $\bgP^1_1\CA_0$ is enough to ensure the existence of (a real coding) the wellfounded ordinal of $\Mii$).  Suppose otherwise; then $\wfo(\Mii) > \alpha$, and then $L_\alpha \in \Mii$.  Then $\fone = \Th(L_\alpha)$ and $\tau$ belongs to $\Mii$, so that $\ftwo = \tau * \fone$ does as well.  As before, this contradicts the assumption that II wins the play; specifically, $\ftwo$ fails to satisfy condition (\ref{VLprojom}).

So $\wfo(\Mii) \leq \alpha$.  We claim $\Mii$ is illfounded.  Otherwise, either $o(\Mii) = \alpha$, in which case we get $\Mii = L_\alpha = \Mi$, in which case (1) holds and I wins; or else $o(\Mii) < \alpha$, so that $\Mii = L_\gamma$ for some $\gamma < \alpha$, so that $\ftwo = \Th(L_\gamma) \in L_\alpha = \Mi$, and again (1) holds, contradicting that $\tau$ is winning for II.

So $\Mii$ is illfounded with $\wfo(\Mii) \leq \alpha$.  Set $\beta = \wfo(\Mii)$.  If there is no model bearing an infinite depth $\Sigma_2$-nesting, then by Lemma~\ref{welchlem} there exists some nonstandard $\Mii$-ordinal $a$, so that, for every nonstandard $\Mii$-ordinal $b$ with $b \leq^{\Mii} a$, we have $T^b_2 \not\subseteq T^\beta_2$.  That is, for all such $b$, there is a $\Sigma_2$ sentence $\sigma$ so that $L_\beta \not\models \sigma$, but $L_b^{\Mii} \models \sigma$.

It is now straightforward to show $\beta,a$ witness the winning condition (2).  Set $a_0 = a$.  Suppose inductively that $a_i$ is a nonstandard $\Mii$-ordinal with $a_i \leq^{\Mii} a$.  Then by choice of $a$, there is some $\Sigma_2$ formula $\sigma$ so that $L_b^{\Mii} \models \sigma$ and $L_\beta \not\models \sigma$; let $\sigma_i$ be the least such under our fixed enumeration of formulae.  If $a_i$ is not limit in $\Mii$, take $a_{i+1}$ to be the greatest limit ordinal of $\Mii$ below $a_i$;  note then $a_{i+1}$ is also nonstandard and below $a$.

Now if $a_i$ is limit in $\Mii$, we have that $\sigma_i$ is of the form $(\exists u)(\forall v)\psi(u,v)$ for some $\Delta_0$ formula $\psi$.  Let $a_{i+1}$ be least so that for some $x \in L_{a_i+1}^{\Mii}$, we have $L^{\Mii}_{a_i} \models (\forall v) \psi(x,v)$.  Then $a_{i+1} <^{\Mii} a_i$, and since $L_\beta \not\models \sigma_i$, we must have that $a_{i+1}$ is nonstandard.  Thus the construction proceeds, and we have that I wins the play $\la \fone, \ftwo \ra$ via condition (2).  So $\tau$ cannot be a winning strategy.
\end{proof}
These claims combine to show that if there is no model with an infinite depth $\Sigma_2$ nesting, then neither player has a winning strategy in the game $G(F)$.  This completes the proof of the theorem.
\end{proof}

We have thus shown that $\Sigma^0_3$ determinacy implies the existence of a model satisfying $\PMI$, and indeed, of some ordinal $\gamma$ so that $L_\gamma \models \PMI$.  The meticulous reader will observe, however, that our proof of determinacy in Section~\ref{sec:det} really only made use of $\nGS$ monotone inductive definitions.  This may at first appear strange, in light of the fact that $\nGS$ is a much smaller class than $\boldsymbol{\Pi}^1_2$.  This situation is clarified somewhat by the following theorem, which shows that if $\gamma$ is minimal with $L_\gamma \models \PMI$, then the $\boldsymbol{\Pi}^1_2$ relations that are correctly computed in $L_\gamma$ are precisely the $\nGS$ relations.

\begin{Theorem} Let $\gamma$ be the least ordinal so that $L_\gamma$ satisfies $\PMI$.  Let $z$ be a real in $L_\gamma$, and suppose $\Phi(u)$ is a $\boldsymbol{\Sigma}^1_2$ formula.  Then there is a $\Game \boldsymbol{\Sigma}^0_3$ relation $\Psi$ so that, for all reals $x$ of $L_\gamma$, we have
   $L_\gamma$ satisfies $\Phi(x)$ if and only if $\Psi(x)$ holds (in $V$, or equivalently, in $L_\gamma$).
\end{Theorem}
\begin{proof} Fix such a formula $\Phi(x)$.  Then there is a recursive tree $T$ on $\om^3$ so that for all $x$, $\Phi(x)$ holds if and only if for some $y$, $T_{\la x,y \ra}$ is wellfounded.  We define a version of the game from Theorem~\ref{IDNfromDet}.  This time, for a fixed real $x$, each player is required to produce their respective $\omega$-models $\Mi, \Mii$ satisfying
 \begin{equation*}\tag{$**$}\label{VLxprojom}
   V=L(x) + KP + \rho_\om = \om.
 \end{equation*}
In addition, $\Mi$ must satisfy the sentence ``$(\exists y)T_{\la x,y \ra}$ is ranked"; whereas $\Mii$ must satisfy its negation.  If a winner has not been decided on the basis of one of these conditions being violated, then Player I wins if either of the conditions (1), (2) from the proof of Theorem~\ref{IDNfromDet} hold.  Let $F_x$ be the set of  $f \in \om^\om$ so that Player I wins the play of the game on $x$, where $\la \fone, \ftwo\ra$, where $\fone(n) = f(2n), \ftwo(n)=f(2n+1)$ for all $n$.  Let $F=\{\la x,f \ra \mid f \in F_x\}$.  Then $F$ is $\Sigma^0_3$; let $\Psi(x)$ be the statement ``I has a winning strategy in the game $G(F_x)$".

Suppose $x \in L_\gamma$ is such that $L_\gamma \models \Phi(x)$.  We claim $\Psi(x)$ holds; that is, Player I has a winning strategy in $G(F_x)$.  Let $y$ be a witness to truth of $\Phi$, and let $\alpha$ be least such that $y \in L_\alpha(x)$ and $L_\alpha(x) \models \KP$.  Then by admissibility, $L_\alpha(x)$ contains a ranking function for $T_{\la x,y \ra}$.  Let $\sigma$ be the strategy for I that always produces the theory of $L_\alpha(x)$.  We claim $\sigma$ is winning for Player I.

Suppose towards a contradiction that $\Mii$ is the model produced by a winning play by II against $\sigma$; we can assume $\Mii \in L_\gamma$.  Then $\Mii$ is an $\omega$-model.  It cannot be wellfounded, since then it would be of the form $L_\beta(x)$ for some $\beta$; but we can't have $\beta \geq \alpha$ (since $\Mii$ cannot contain $y$, or else it would have a ranking function for $T_{\la x,y\ra}$), nor can $\beta < \alpha$ hold (since then (1) is satisfied, and I wins the play).  So $\Mii$ is illfounded, say with $\wfo(\Mii) = \beta$; by a similar argument, $\beta < \alpha$.  Now since I does not win the play, the condition (2) fails, so there must be some infinite depth $\Sigma_2$-nesting based on $\Mii$, by Lemma~\ref{welchlem}.  But this contradicts the fact that the model $\Mii$ belongs to $L_\gamma$, by minimality of $\gamma$ and Theorem~\ref{PMIfromIDN}.

Conversely, suppose $\Phi(x)$ fails in $L_\gamma$.  Suppose towards a contradiction that Player I wins the game $G(F_x)$ (in $V$); then by Theorem~\ref{detfrommi}, there is such a strategy $\sigma \in L_\gamma$.  Let $\Mii$ be the least level of $L(x)$ containing $\sigma$.  Note that $\Mii \models ($\ref{VLxprojom}$)+``(\forall y) T_{\la x,y \ra} $ is not ranked".  By the argument in the proof of Theorem~\ref{IDNfromDet}, we obtain failure of both (1) and (2), so that II wins the play, a contradiction to $\sigma$ being a winning strategy.
\end{proof}

\bibliographystyle{asl}
\bibliography{thesbib}
\end{document}